\documentclass{article}

\usepackage {amssymb}

\usepackage {amsfonts}

\usepackage {amsmath}

\usepackage {amsthm}

\usepackage{graphics}

\usepackage {framed}

\usepackage {amsxtra}

\usepackage {enumerate}

\usepackage[margin=1in]{geometry}

\theoremstyle{plain}

\newtheorem{thm}{Theorem}

\newtheorem{cor}[thm]{Corollary}

\newtheorem{lem}[thm]{Lemma}

\newtheorem{prop}[thm]{Proposition}

\theoremstyle{definition}

\theoremstyle{remark}

\begin{document}

\title{On powers of the Diophantine function $\star:x\mapsto x(x+1)$}

\author {Maxwell Lippmann\footnote{max.lippmann@yahoo.com\ \ Roslyn High School\ \ Roslyn, NY 11576 } \and Donald Silberger\footnote{SilbergD@newpaltz.edu or, preferably, DonaldSilberger@gmail.com\ \ State University of New York\ \ New Paltz, NY 12561}}

\maketitle

\[\mbox{\sf Dedicated to Jan Mycielski,\,\, 1932 February 07 - 2025 January 23}\]

\begin{abstract}

We treat the functions $\star^k:{\mathbf N}\rightarrow{\mathbf N}$ where $\star:x\mapsto \star x := x(x+1)$. The set $\{\star^k x+1: \{x,k+1\}\subset{\mathbf N}\}$ is pairwise coprime; so, the set ${\mathbf P}$ of primes is infinite. Our Theorem \ref{Embed} resorts to the {\it mother sequence}, {\sc m},  that is obtained by factoring the infinite integer-sequence $2,3,4,5,\ldots$ into prime powers. 

For each $x\ge1$ we define the {\it gross $x$-sequence}, $\gamma_\star(x) :=\langle x+1;\star x+1;\star^2x+1;\star^3x+1;\ldots\rangle$, and also the {\it star sequence}, $x^\star$, obtained by factoring the terms of $\gamma_\star(x)$ into prime powers.  It turns out that $\gamma_\star(1)$ is Sylvester's sequence, A00058 in the {\sl On-Line Encyclopedia of Integer Sequences}, OEIS, and that $\gamma_\star(2)$ is the sequence A082732 in the OEIS.

{\bf Theorem} \ref{NotAll}. For every integer $x\ge1$ there is a prime $p(x)$ that divides no member of $\{\star^kx+1: k\ge0\}$.

{\bf Theorem}  \ref{Embed}. For each sequence $\eta$ of powers of primes there  are infinitely many subsequences ${\sf c}_j$ of {\sc m} such that numerically $\eta = {\sf c}_j$ but where the term-set family in {\sc m} of those ${\sf c}_j$ is formally pairwise disjoint.

{\bf Theorem} \ref{Recip}.  $1/x = \sum_{k=0}^{n-1} 1/(\star^kx+1)+1/(\star^nx) = \sum_{k=0}^\infty 1/(\star^kx+1)$ for all $\{x,n\}\subseteq{\mathbf N}$.

{\bf Theorem} \ref{Diverge}. For every $x\in{\mathbf N}$, when $x^\star := \langle x_j\rangle_{j=0}^\infty$ then $\sum_{j=0}^\infty 1/x_j  =  \infty$.
\end{abstract}

\section{Introduction}

Our initial interest in the elementary function $\star:{\mathbf N}\rightarrow{\mathbf N}$ defined by $\star:x\mapsto \star x := x(x+1)$, where ${\mathbf N} := \{1,2,3,\ldots\}$, derives from its appearance in what we have called in \cite{Silbergers,Silberger} the Vital Identity

\[ \frac{1}{z} = \frac{1}{z+1} + \frac{1}{z(z+1)}\quad\mbox{This identity holds for every complex number}\quad z\notin \{-1,0\}. \]

\noindent In \cite{Silbergers, Silberger} the Vital Identity has facilitated proofs about segments of the harmonic series and more generally about sums of the reciprocals of positive integers. We define $\star^0$ to be the identity function and $\star^{k+1} := \star^k\circ\star$.

We recall Theorem 1 in \cite{Silberger} where the fact that two consecutive positive integers are necessarily coprime entails that for every pair $k\ge0$ and $x\ge1$ of integers the integer $\star^k x+1$ has at least $k$ distinct prime factors. This fact immediately yields Euclid's theorem that there exist infinitely many prime integers.\footnote{This inescapable observation arising from Sylvester's sequence has been noted by other people; e.g., R. W. K. Odoni mentions it  in the second sentence of his 1984 paper \cite{Odoni}.}
\vspace{.3em}

We are interested  for each $x\in{\mathbf N}$ in two infinite sequences of integers. One is the {\it star sequence} $x^\star$; its terms are  powers of primes. Each such  $x^\star$ is uniquely determined and is created thus: First we produce the other sequence of interest to us; this is the {\it gross $x$-sequence}  \[ \gamma_\star(x) := \langle x+1;\star x+1;\star^2 x+1;\star^3x+1;\ldots;\star^kx+1;\ldots\rangle = \langle \star^kx+1\rangle_{k=0}^\infty\, . \]The terms of $\gamma_\star(x)$ are separated by semicolons. It is easy to see, but it also follows from Theorem 3 of \cite{Silberger}, that $\gamma_\star(x)$ is injective, that it very rapidly increases and that the set of its terms is pairwise coprime. 

We manufacture $x^\star$ from  $\gamma_\star(x)$ by factoring each term of $\gamma_\star(x)$ into a product of the powers of its primes and listing them  separated by commas between successive semicolons in ascending order of prime sizes. Change the semicolons into commas. The integers separated by commas are the terms of $x^\star := \langle x_i\rangle_{i=1}^\infty$. \vspace{.3em}

\noindent{\bf Example One.} The gross $1$-sequence $\gamma_\star(1)$ begins with  $\langle 2;3;7;43;1807;3263443;\ldots\rangle$ and turns out to be  Sylvester's sequence, A00058 in the {\sl On-Line Encyclopedia of Integer Sequences - OEIS}. We compute that $1^\star = \langle 2,3,7,43,13,139,3263443,\,\ldots\rangle$; i.e., $1_0=2,\, 1_1 = 3;\,1_2=7,\, 1_3=43,\, 1_4=13,\, 1_5=139,\, 1_6=3263443$. 
The OEIS says that it is an open question whether every term in $\langle w_k \rangle_{k=1}^\infty$  is square-free.\vspace{.3em}

\noindent{\bf Example Two.} We start the gross $2$-sequence $\gamma_\star(2) = \langle 3; 7; 43; 1807; 3263443; \dots \rangle$ and then we compute $2^\star = \langle 3,\,  7,\, 43,\, 13,\, 139,\, 3263443,\,\ldots\rangle$. Our $\gamma_\star(2)$ occurs as A082732 in the {\sl OEIS.} \vspace{.3em}

\noindent{\bf Example Three.} The gross $3$-sequence begins with $\gamma_\star(3) = \langle 4; 13; 157; 24493; 599882557;\ldots\rangle$. We compute $3^\star = \langle 2^2,\,13,\, 157,\, 24493,\, 67,\, 277,\,32323,\,\ldots\rangle$\vspace{.3em}

Notice that $\gamma_\star(2)$ is a proper suffix of $\gamma_\star(1)$. This implies that $2^\star$ is a proper suffix of $1^\star$. These observations together with Theorem 3 in \cite{Silberger} give us   

\begin{prop}\label{Prefix} If $y\not= x$ and $y+1$  is a term in the sequence $\gamma_\star(x)$, then $y^\star$ is a proper suffix of $x^\star$, and no prime that is a factor of a term of the prefix of $x^\star$ that is complementary to $y^\star$ is a factor of any term of $y^\star$.\end{prop} 

Proposition \ref{Prefix} informs us that there are integers $x\ge1$ for which there are primes that are factors of no member of $\{\star^k x+1: k\ge0\}$.  Theorem \ref{NotAll}, below, extends this fact.

For $x\in{\mathbf N}$ we define the sequence $\eta(x) := \langle w_{x,j}\rangle_{j=1}^\infty$  recursively by \[ w_{x,k+1} = 1 + \prod_{j=0}^k w_{x,j}\quad\mbox{and}\quad w_{x,0} := x. \] According to Odoni  \cite{Odoni}, whose paper's title contains it, $\eta(1)$ is Sylvester's sequence, A00058 in the OEIS.

\begin{lem}\label{Sylvester} Let $x\in{\mathbf N}$. The sequence $\eta(x) := \langle w_{x,j} \rangle_{j=1}^\infty$ is identical to the sequence $\gamma_\star(x) := \langle \star^kx+1 \rangle_{k=0}^\infty$. \end{lem}

\begin{proof}  We show that a variant of the recursion defining $\eta(x)$ produces $\gamma_\star(x)$.  Observe that $\star^0x+1 = x+1 = w_{x,1}$. Suppose for $k>0$ that $w_{x,k} = \star^{k-1}x+1$. Then $\star^kx +1 := (\star^{k-1}x)(\star^{k-1}x+1) +1 = (w_{x,k}-1)w_{x,k} +1 = (w_{x,1}w_{x,2}\cdots w_{x,k-1})w_{x,k}+1 = 1+ w_{x,1}w_{x,2}\cdots w_{x,k} = w_{x,k+1}$. The lemma follows by induction. \end{proof}

\begin{thm}\label{NotAll} For every integer $x\ge1$ there is a prime $p(x)$ that divides no member of $\{\star^k x+1: k\ge0\}$.  \end{thm}

\begin{proof} In his paper's  first paragraph R. W. K. Odoni \cite{Odoni} shows that, if $p$ is a prime that divides a term of $\eta(1)$, then either $p=3$ or $p\equiv1\mod 6$.  Therefore $11$ divides no term $w_k$ of $\eta(1)$, and hence by Lemma \ref{Sylvester}  we have that $p(1) := 11$ is a factor of no term of the sequence $\gamma_\star(1) := \langle \star^j1+1\rangle_{j=0}^\infty$. 

Let $1 < x \in{\mathbf N}$, and let $p(x)$ be any prime factor of $x$. Observe that $x$ is coprime to every term of the sequence $\gamma_\star(x)$ and therefore that $p(x)$ divides no term of $\gamma_\star(x) := \langle \star^jx+1 \rangle_{j=0}^\infty$.  \end{proof}

\section{Parallel embeddings}

Our $x^\star$ are members of the class of sequences $\eta$ whose terms are powers of prime integers. Each of these $\eta$  occurs in an  infinite assemblage of subsequences of the {\it mother sequence} {\sc m} $:= \langle m_0,m_1,m_2.\ldots\rangle = \langle m_i\rangle_{i=0}^\infty$ in whose womb each $\eta$ resides along with infinitely many pairwise formally disjoint identical siblings ${\sf c}_j = \eta$. 

The sequence {\sc m} is created by listing the natural numbers greater than $1$ in their normal ascending order but then replacing each of those integers with its normal-order factorization into powers of primes separated by commas. The following finite prefix of the infinite sequence {\sc m} ought to eliminate  misunderstandings:
\[\mbox{\sc m} = \langle 2,3,2^2,5,2,3,7,2^3,3^2,2,5,11,2^2,3,13,2,7,3,5,2^4,17,2,3^2,19, 2^2,5,3,7,2,11,23,2^3,3,5^2,2,13,3^3,\ldots\rangle\]

Notice that $2 = m_0 = m_4 = m_9 = m_{15} = m_{21} = \cdots$ and that $5 = m_3 = m_{10} = m_{18} = \cdots$ That is, not only does {\sc m} fail to be injective, but every prime power occurs infinitely often as a term in {\sc m}. We call the symbols $m_i$ and $m_j$ {\it formally distinct} as terms of {\sc m} if and only if $i\not=j$. We emphasize that $m_i$ and $m_j$ can be formally distinct even in the event that a numerical equality $m_i = m_j$ obtains. For instance, we say that $m_0$ ``formally distinct'' from $m_4$ although both $m_0$ and $m_4$ are numerically identical; i.e., $m_0=2=m_4$. 

We call a collection ${\bf H}$ of subsequences of {\sc m} {\it pairwise formally disjoint} - aka {\it PFD} - if and only if, for every pair $\{{\sf c,d}\}\subset{\bf H}$, whenever $\langle m_i,m_j\rangle\in{\cal T}({\sf c})\times{\cal T}({\sf d})$ then $i\not= j$, where ${\cal T}({\sf c})$ is the term set of ${\sf c}$.

By a {\it parallel embedding} of a sequence $\eta$ of powers of primes we mean an infinite PFD collection of subsequences of {\sc m} each of which is term by term a numerical copy of $\eta$. 

\begin{thm}\label{Embed} Every sequence $\eta$ of powers of primes has a parallel embedding.
 \end{thm}
 
\begin{proof} We leave to the reader the very easy proof that every finite $\eta$ has a parallel embedding, and prove the theorem for infinite $\eta$ using the paradigm argument pertaining to $\eta := x^\star$ for an arbitrary $x\in{\mathbf N}$. Recall from Section 1 the terminology $x^\star := \langle x_0,x_1,x_2,\ldots\rangle$. We will recursively develop the promised parallel embedding ${\bf S}(x) := \{ {\sf c}_0, {\sf c}_1, {\sf c}_2,\ldots\}$ consisting of subsequences ${\sf c}_j$ of {\sc m} that term-by-term satisfy the numerical identity ${\sf c}_j = x^\star$ for all $j\ge0$. In order for our procedure to build all of the ${\sf c}_j\in{\bf S}(x)$ and not merely a solitary one such infinite sequence, we shall resort to the trick used to prove Theorem 2 in \cite{Silberger} and reminiscent of the corner-slicing gimmick employed by Georg Cantor \cite{Cantor} in his proof that the set ${\mathbf Q}^+$ is countable.

We intend each ${\sf c}_i$ to be a sequence ${\sf c}_i := \langle m_{i_0}, m_{i_1}, m_{i_2},\ldots, m_{i_j},\ldots \rangle$ of powers of primes for which $m_{i_j} = x_j$ is a numerical fact for every $j\ge0$. Let $0_0$ be the least index $t$ for which $m_t= x_0$. The term $m_{0_0}$ was just now ``captured from {\sc m}.'' For our purposes in building the set ${\bf S}(x)$ of subsequences of {\sc m}, our designation of $m_{0_0}$ as ``captured'' in effect removes $m_{0_0}$ as an available term from the mother sequence {\sc m}.

Next, let $0_1$ be the smallest integer $t$ such that $m_t$ is not yet captured and also such that $m_t = x_1$. We have thus far chosen the first two terms, $m_{0_0}$ and $m_{0_1}$, of the subsequence ${\sf c}_0$. Our third choice will be the first term $m_{1_0}$ of the subsequence ${\sf c}_1$; here $1_0$ is the smallest integer $t$ for which both $m_t$ is uncaptured and $m_t = x_0$. Three subsequential terms have been designated and thus captured. The order in which we select (and capture) terms of {\sc m} for subsequences destined to comprise ${\bf S}(x)$ is easily recognized from this pattern:

$m_{0_0}\rightarrow m_{0_1} \rightarrow m_{1_0} \rightarrow m_{0_2} \rightarrow m_{1_1} \rightarrow m_{2_0} \rightarrow m_{0_3} \rightarrow m_{1_2} \rightarrow m_{2_1} \rightarrow m_{3_0} \rightarrow m_{0_4} \rightarrow m_{1_3} \rightarrow m_{2_2} \rightarrow \cdots$

As a paradigm step in our subsequence term-choosing procedure, we now choose the second term of the subsequence ${\sf c}_3$. Let $3_1$ be the smallest integer $t$ for which $m_t=x_1$ but also for which $m_t$ is not yet captured. We capture $m_{3_1}$ from {\sc m} to serve as the the second term, $m_{3_1}$, in the subsequence ${\sf c}_3$.

It is clear for every $\langle u,v\rangle\in\omega\times\omega$ that\footnote{Following the custom of set theoreticians we define  $\omega := \{0\}\cup{\mathbf N}$.}, after a finite number of procedural steps of the sort described in the preceeding three paragraphs, the $(v+1)^{\rm st}$ term, $m_{u_v}$, of the subsequence ${\sf c}_u$ will be captured and put into its place in ${\sf c}_u$ and that $m_{u_v} = x_v$. Moreover, for each pair $p \not= q$ of nonnegative integers, our refusal to re-use any previously captured term of {\sc m}, that term's having already been assigned to be a term of some ${\sf c}_r$, guarantees that the set of terms of the subsequence ${\sf c}_p$ is formally disjoint from the set of terms of ${\sf c}_q$.\end{proof}

\begin{cor}  Every copy of {\sc m} has a  parallel embedding the family of whose term sets is a partition of {\sc m}. \end{cor}

\section{Sequence sums}

In \cite{Silbergers,Silberger} the function $\sigma:{\cal P}({\mathbf N})\setminus\{\emptyset\}\rightarrow {\mathbf Q}^+\cup\{\infty\}$ is defined by $\sigma: X\mapsto \sigma X := \sum_{y\in X} 1/y$, where ${\cal P}({\mathbf N})$ is the family of sets of positive integers. 

Define ${\rm T}(x)$ to be the set of terms $x_i$ of $x^*$. Remember that each such  $x_i$ is the power of a prime.

Our next result generalizes a fact about Sylvester's sequence, i.e., the case $x=1$.

\begin{thm}\label{Recip} Let $x\ge1$. Then $1/x = \sigma\bigcup_{k=0}^n\{\star^k x+1\} + 1/\star^nx$ for all $n\in{\mathbf N}$ and $1/x = \sigma\bigcup_{k=0}^\infty \{\star^k x+1\}$.
 \end{thm}

 \begin{proof} Since $\star^kx\notin\{-1,0\}$ when $k\ge0$, by the Vital Identity we have that  \[ \frac{1}{x} = \frac{1}{x+1}
+ \frac{1}{\star x} =: \frac{1}{\star^0x+1} + \frac{1}{\star^1x} = \frac{1}{\star^0x+1}+ \frac{1}{\star^1x+1} + \frac{1}{\star^2x} =  \frac{1}{\star^0x+1} +  \frac{1}{\star^1x+1} + \frac{1}{\star^2x+1} + \frac{1}{\star^3x} = \]  
\[ \frac{1}{\star^0x+1} + \frac{1}{\star^1x+1} + \frac{1}{\star^2x+1} + \frac{1}{\star^3x+1}+\frac{1}{\star^4x} =  \,\,\cdots\,\, = \sum_{k=0}^{n-1} \frac{1}{\star^kx+1} + \frac{1}{\star^nx}.\] The first claim is proved. The second follows from it together with the fact that $\lim_{n\rightarrow\infty} 1/\star^nx = 0$. \end{proof}

In the light of Theorem \ref{Recip}, of the paucity of factors that occur in each of the first few terms of the sequences $\gamma_\star(x)$ at which we have looked, and also of the swiftly burgeoning sizes of the primes encountered there, one might wager that $\sigma{\rm T}(x)<\infty$. However, there are reasons to hedge this bet.

One such reason is that $1/xy < (x+y)/xy = 1/x + 1/y$ with $x>1$ and $y>1$. So, if  factors-rich terms of the sort $\star^kx +1 = p_{i_1}^{e_1}p_{i_2}^{e_2}\ldots p_{i_v}^{e_v}$ occur frequently in $\gamma_\star(x)$, then such numbers as $1/(\star^kx+1)$ that are summands in the sum treated in Theorem \ref{Recip} get replaced by significantly larger numbers  $1/p_{i_1}^{e_1}+\ldots +1/p_{i_v}^{e_v}$ that are   partial sums in $\sigma{\rm T}(x)$. This makes the convergence of $\sigma{\rm T}(x)$ begin to seem questionable. For, if among the terms of $\gamma_\star(x)$ there is an increasing manifestation of factors-rich terms as $k$ increases, then $\sigma\bigcup\{\star^kx+1: k\ge0\} \ll \sigma{\rm T}(x)$. Indeed,  $\sigma {\rm T}(x) = \infty$ appears possible. 

On the other hand, since huge primes do occur early as terms in some sequences $x^\star$, possibly prime-power terms of  $x^\star$ tend to increase in size swiftly enough to entail rapidly diminishing summands in $\sigma{\rm T}(x)$ so as to support  convergence,  $\sigma{\rm T}(x)<\infty$, after all. 

So, what actually happens? 

\begin{thm}\label{Diverge} $\sigma{\rm T}(x) = \infty$ for every integer $x\ge1$. \end{thm} 

\begin{proof} We rely upon the Second Theorem of Franz Mertens \cite{Mertens} and upon Theorem 2 of R. W. K. Odoni  \cite{Odoni} together with Section 8 of Odoni's paper where he observes that the provisions of his Theorem 2 that are met by $\gamma_\star(1)$ hold also for the sequences $\gamma_\star(x)$ with $x>1$.  

The terms of $x^\star := \langle x_j \rangle_{j=0}^\infty$ are powers of primes. We write $x_j^-$ to designate the prime factor of $x_j$.

The Mertens theorem estimates the size of the sum of the reciprocals of all primes less or equal to $n$ by \[\lim_{n\rightarrow\infty}\Big(\sum_{p\le n} \frac{1}{p} - \log\log n - M\Big) = 0 \quad\mbox{where}\quad M := 0.26149721\ldots \quad\mbox{is the Meissel-Mertens constant.}\]

The Odoni theorem estimates the number $\underline{\bf P}_1(n)$ of distinct prime factors less than or equal to $n$ of terms, $\star^k1+1$ of the sequence $\gamma_\star(1)$, by asserting that  \[\underline{\bf P}_1(n) = O\Big(\frac{n}{(\log n)( \log\log\log n)}\Big)\quad\mbox{as}\quad n\rightarrow\infty.\]

Section 8 of the Odoni paper claims that the foregoing assertion about $\gamma_\star(1)$ applies to $\gamma_\star(x)$ for all $x\ge1$ and consequently extends our statements here about the sequence $1^\star$ to the sequences $x^\star$ as well.

Since the Prime Number Theorem estimates $n/\log n$ to be the number $\pi(n)$ of primes no greater than $n$ as $n\rightarrow\infty$,  the Odoni paper implies that the proportion of $\pi(n)$ of those such primes that occur as factors of terms of $\gamma_\star(x)$ is asymptotic to $1/(\log\log\log n)$. Using also the Merten theorem we infer that  \[ \sum_{p\in\underline{\bf P}_x(n)} \frac{1}{p} \approx \frac{\log\log n +M}{\log\log\log n}\,\,\mbox{whence}\,\,\sum_{j=0}^\infty\frac{1}{x_j^-} =   \sum_{p\in\underline{\bf P}_x}\frac{1}{p} = \infty\,\,\mbox{since}\quad\lim_{n\rightarrow\infty}\frac{\log\log n +M}{\log\log\log n} = \infty,\,\,\mbox{where}   \]  $\underline{\mathbf P}_x$ is the set of all prime factors $p$ of terms of $\gamma_\star(x)$ and $\underline{\mathbf P}_x(n)$ is the set of such $p$ no larger than $n$. Thus although we have obtained the tightening \[ \sum_{j=0}^\infty \frac{1}{x^-} = \infty \] of Leonhard Euler's 1937 theorem \cite{Euler1} that $\sum_{p\in {\mathbf P}} 1/p = \infty$, we cannot yet infer as ultimately desired that
 \[ \sigma{\rm T}(x) := \sum_{j=0}^\infty \frac{1}{x_j} = \infty. \]

We must dismiss the chance that $x_j = x_j^-$ happens too seldom to compel $\sum_{j\in\omega} 1/x_j$ to diverge. That this might  threaten our theorem is suggested by Euler's Basel problem solution \cite{Euler2}, $\sum_{n\in{\mathbf N}} 1/n^2= \pi^2/6$.

We will  use probability theory to make an end run around this ostensible hazard we have noted.\vspace{.3em}

It is evident that $\sigma{\rm T}(x) = \infty$ if and only if, for every probability $r<1$ and every real number $\ell$, there is an index $b$ such that the probability is greater than or equal to $r$  that $\sum_{j=b}^\infty 1/x_j > \ell$.

Since no prime divides more than one term $x_j$ of the sequence $x^\star$, no prime power can recur with an increased exponent.  Moreover, among all multiples $z\in{\mathbf N}$ of a prime $p$, the probability that $p^2|z$ is equal to $1/p$ whence the probability that $\neg (p^2|z)$ is equal to $1-1/p = (p-1)/p$. 

For each $v\in{\mathbf N}$ there is an index $b$ such that the smallest prime dividing any term in the suffix subsequence $\langle x_{b+i} \rangle_{i=0}^\infty$ of $x^\star$ exceeds $v$. So, for $r<1$ a probability, there exists $b\in{\mathbf N}$ with $r<(x_j^--1)/x_j^- <1$ for every $j\ge b$. Hence, for some $m>b$ we have that the fraction, $f(b,m)$ of the number of terms of the sort  $x_j=x_j^-$ divided by the number $m-b+1$, is greater than $r$. Thus, $\liminf_{m\rightarrow\infty} f(b,m) \ge r$. It follows, for an arbitrary real number $\ell$, that there exists $b$ for which the probability is greater than or equal to $r$ that $\sum_{j=b}^\infty 1/x_j > \ell$. We conclude that $\sigma {\rm T}(x) = \infty$. \end{proof}

The proof of Theorem \ref{Diverge} gives us the following fact gratis.

\begin{cor}\label{Pi} For all $x\in{\mathbf N}$ the approximation $\sum_{p\in\underline{\bf P}_x(n)} \frac{1}{p} \approx \pi(\log\log n)$ obtains as $n\rightarrow\infty$.   \end{cor}
\vspace{.5em}

Seeking a boundary between convergence and divergence we pose two questions pertaining to each $x\in{\mathbf N}$: 
\vspace{.1em} 

\noindent{\bf One.} If  ${\sf b}_{x,k}$ is the largest prime with ${\sf b}_{x,k}^{e_{x,k}}\|\star^k x+1$ for $e_{x,k}>0$, then  $\sigma\{{\sf b}_{x,k}^{e_{x,k}}: k\ge0\} < \infty$?\vspace{.5em}

\noindent{\bf  Two.} If ${\sf d}_{x,k}$ is the smallest prime with ${\sf d}_{x,k}^{f_{x,k}}\|\star^kx+1$ for $f_{x,k}>0$, is $\sigma\{
{\sf d}_{x,k}^{f_{x,k}}:k\ge0\}<\infty$?\vspace{1em}

Googling Sylvester's sequence led us to Wikipedia which offers hundreds of related references, and it was Wikipedia that apprised us of the valuable {\sl On-Line Encyclopedia of Integer Sequences} where we found that our $\gamma_\star(1)$ had been discovered in another form  by James Joseph Sylvester in 1880. Indeed, Euclid may have known about that sequence two and one-third millennia ago. 

Both $\gamma_\star(1)$ and $\gamma_\star(2)$ (respectively A00058 andA082732 in the {\sl OEIS}) have number-theoretic significance. Each term of  $\gamma_\star(1)$ is $1$ plus the product of the terms that precede it and each term of $\gamma_\star(2)$ is $1$ plus the least common multiple of the terms preceding it. The OEIS says also that, starting with the fifth term of the sequence $\gamma_\star(2)$, the ultimate two-digit suffixes of its terms alternate perpetually between 57 and 93.\vspace{.5em}

Are there  number-theoretic interpretations of $\gamma_\star(x)$ for other values of $x$ besides $x=1$ and $x=2$? 
For $x$ a prime power, what patterns emerge for the $\gamma_\star(x)$ and the $x^\star$?  What happens when $x$ is square-free?\vspace{.5em}

The set $\{\star^k(1+\prod{\mathbf P}_{\rm known}): k\ge0\}$ supplements $1+ \prod{\mathbf P}_{\rm known}$ as a source for undiscovered primes.

\vspace{1em}
\noindent{\bf Acknowledgments.}  Kira Adaricheva, Evan O'Dorney and Allan J. Silberger offered us helpful counsel.

\vspace{2em}

\noindent{\bf 2020 Mathematics Subject Classification:} \, 01A55, 01A60, 11A25, 11N05, 11N25

\end{document}